\documentclass[a4paper,12pt]{article}

\usepackage{amsthm}
\usepackage{amsmath,amssymb,latexsym,amsfonts,mathrsfs}
\usepackage[dvipdfmx]{graphicx}

\usepackage{color}

\usepackage{fancyhdr}
\usepackage[top=30truemm, bottom=30truemm, left=25truemm, right=25truemm]{geometry}

\newcommand{\ep}{\varepsilon}
\newcommand{\nn}{\nonumber}

\newcommand{\SCR}[1]{{\mathscr #1}}

\newcommand{\CAL}[1]{{\cal #1}
}

\newcommand{\J}[1]{\left\langle #1 \right\rangle}
\newcommand{\D}[1]{{\mathscr D}( #1 )}

\theoremstyle{definition}
\newtheorem{Thm}{{\bf Theorem}}[section]

\newtheorem{Lem}[Thm]{{\bf Lemma}}

\newtheorem{Cor}[Thm]{{\bf Corollary}}
\newtheorem{Ass}[Thm]{{\bf Assumption}}

\newtheorem{Rem}[Thm]{{\bf Remark}}

\newcounter{Exami}

\newcommand{\Proof}[2][Proof]{
\begin{proof}[{\bf #1}]
#2
\end{proof}
}

\begin{document}

\begin{flushleft}
{\bf \Large High-Energy Smoothing Estimates for Selfadjoint Operators
 } \\ \vspace{0.3cm} 
by {\bf \large  Masaki Kawamoto 
 $^{1}$} \\ 
 Department of Mathematics, Faculty of Science, Tokyo University
 of Science, 1-3, Kagurazaka, Shinjuku-ku, Tokyo 162-8601, Japan. \\
 Email: mkawa@rs.tus.ac.jp
\end{flushleft}

\begin{center}
\begin{minipage}[c]{400pt}
{\bf Abstract}. {\small 
We prove the limiting absorption principle on the non-compact interval $I$, on which the uniformly positive Mourre estimate holds. We reveal that such a result yields so-called smoothing estimates.
}
\end{minipage}
\end{center}

\begin{flushleft}
{\bf Keywords}; Mourre theory; Kato's smoothing; Limiting absorption principle; scattering theory.
\end{flushleft}
\section{Introduction}
In this paper, we show the limiting absorption principle (LAP) for generalized selfadjoint operators $H$ by using Mourre's theory. In the paper by Mourre \cite{Mo}, the energy localized LAP was deduced for $H$. Hence in this paper, we deduce the global in the high-energy LAP and a smoothing-type LAP by modifying Mourre's approach. The key to doing so is the so-called positive commutator $A$ (called the conjugate operator), such that $f(H) i[H,A] f(H) \geq c f(H) ^2$ with $c>0$ and energy cut-off $f \in C_0^{\infty} ({\bf R})$. In the background of Mourre theory, a pair of Schr\"{o}dinger operators $H_S =  p^2 $ are considered, along with the generator of dilation group $A_S= x \cdot p + p \cdot x$, where $x = (x_1,...,x_n) \in {\bf R}^n$ and $p = -i \nabla _x$. Formally, it then holds that $f(H_S) i[H_S,A_S] f(H_S) =  f(H_S) ( 4p^2 -2 \nabla V)f(H_S) $ on $L^2({\bf R}^n)$. By selecting $f $ such that $\{ 0 \} \cup \{ \sigma_{\mathrm{pp}} (H) \} \notin \mathrm{supp} (f) $, we have $f(H_S) 4p^2 f(H_S) \geq 2c_0 f(H_S)^2 $, $c_0 >0$, and by selecting a very narrow $f$, $f(H_S) |2\nabla V| f(H_S) \leq c_0 f(H_S)^2 $ holds if $V$ and $f$ satisfy the suitable conditions. Consequently, we can deduce that $f(H_S)i[H_S,A_S] f(H_S) $ is positive with a suitable energy cut-off. The condition for the support of $f$ is mainly used to deduce the smallness of $Vf(H_S)$ and to deduce the positiveness and boundedness of the commutator $f(H_S) i[H_S,A_S] f(H_S)$. On the other hand, the Yokoyama-type conjugate operator $A_Y := x \cdot p (1+p^2)^{-1} + (1+p^2)^{-1} p \cdot x$ (see Yokoyama \cite{Yo}) gives $i[H_S ,A_Y] = 4 p^2 \J{p}^{-2} + \J{p}^{-2} \times (\mbox{bounded operators})$ if $V$ is smooth and satisfies a suitable decaying condition, where $\J{\cdot} = (1+ \cdot ^2)^{1/2}$. For the high-energy cut-off $f(H_S)$ (i.e., $p^2 \gg 1$ holds on $f(H_S)$), it holds that $f(H_S) i[H_S ,A_Y]  f(H_S) \geq 3 f(H_S)^2$ even if the support of $f$ is not narrow. Moreover $i[H_S, A_Y]$ can be extended to a bounded operator without the energy cut-off. Additionally, in many physical situations the Hamiltonian $H$ has a conjugate operator $A$, which gives the positive and bounded commutator for the high-energy cut-off. Judging from these, we expect that we can deduce the limiting absorption principle from Mourre’s theory with a pair of $H$, $A$ and $\varphi$ such that $i[H,A]$ is positive and bounded with a high-energy cut-off $\varphi(H)$, with $\varphi \in C^{\infty} ((-\infty ,-R] \cup [R, \infty))$ (not $C_0^{\infty} ({\bf R})$), where $R$ is a given large constant. 

Consider the Hamiltonian $H$, where $H$ is a selfadjoint operator acting on Hilbert space $\SCR{H}$. The norm on $\SCR{H}$ and  the operator norm on $\SCR{H}$ are denoted by the same notation, $\| \cdot \|$, and $(\cdot , \cdot)$ denotes the inner product of $\SCR{H}$. Let $A$ be a selfadjoint operator and suppose $D= \D{H} \cap \D{A} \subset \SCR{H}$ is dense. We define the form $q_{H,A} (\cdot, \cdot)$ on $D$ as $q_{H,A} (u,v) := i(Au,Hv) -i(Hu,Av)$ for $u,v \in D$. 
Then, we say that $i[H,A]$ can be extended to a bounded operator from $\SCR{H}$ to $\SCR{H}$ if there exists a bounded selfadjoint operator $T$ such that the closure of $q_{H,A} (\cdot , \cdot)$, $\tilde{q}_{H,A}(\cdot,\cdot)$, satisfies $\tilde{q}_{H,A} (u,v) = (Tu,v)$, $u,v \in \SCR{H}$, and denote this by $T = i[H,A]^0$. We now introduce the first assumption on $H$, $A$ and energy cut-off $\varphi$ in order to derive the high-energy LAP:
\begin{Ass} \label{A1}
Let $H$ and ${A}$ be selfadjoint operators satisfying the following statements: (a), (b), (c), (d), (e), and (f) (all conditions stated in Mourre \cite{Mo} are fulfilled, regarding which see also \cite{CFKS}): 
\\ 
(a). $\overline{H|_{\D{H} \cap \D{A}}} =H$. \\ 
(b). $e^{itA} \D{H} \subset \D{H}$, $|t| <1$ and for $u \in \D{H}$, $\sup_{|t| <1} \| He^{itA} u\| < \infty$. \\ 
(c). The commutator $i[H,A]^0$ can be defined in the form and can be extended to a bounded operator from $\SCR{H}$ to $\SCR{H}$. \\ 
(d). There exists a constant $c_2 >0$ such that for all $u,v \in \D{A}$, 
\begin{align*}
\left| 
(Au, i[H,A]^0v) - (i[H,A]^0u, Av)
\right|  \leq c_2 \| u \| \| v \|.
\end{align*}
(e). Let $R>0$ be a given constant and $\varphi =\varphi _R \in C^{\infty} ({\bf R})$ with $0 \leq \varphi \leq 1$ is smooth cut-off function so that $\varphi (s) = 0$ if $ |s| \leq R $ and $\varphi (s) = 1$ if $|s| \geq 2R$. Then there exists $0< c_0 \leq c_1$, such that 
\begin{align*}
c_0 \varphi (H) ^2 \leq \varphi (H) i[H,A] ^0 \varphi (H) \leq c_1 \varphi(H) ^2 \leq c_1.
\end{align*}
Moreover, $\varphi ' (s)  >0$ for all $s$ in $R< |s| < 2 R$.  \\
(f). $i[H,i[H,A]^0]^0$ and $i[i[H,A]^0,H]^0$ can be extended to the bounded operators.   
\end{Ass}
Under this assumption, we have the Mourre estimates for the high-energy case, which provide the following nonexistence-of-eigenvalues property and limiting absorption principle for $H$: 
\begin{Thm}\label{T1}
Under Assumption \ref{A1}, $H$ has no eigenvalues on $(-\infty, -3R] \cup [3R, \infty)$. 
\end{Thm}
\begin{Thm}\label{T2}
For any fixed $0< \hat{\ep} <1$, we select $\tilde{R}$ arbitrarily in $(0, R^{1-\hat{\ep}}]$. Under Assumption \ref{A1}, for large $R$, there exists a constant $\delta = \delta _R >1$ with $\delta _R \to 1$ (as $R \to \infty$) such that the limiting absorption principle 
\begin{align} \nn
& \sup_{| \lambda| \geq 3R, \ \mu >0} \left\| \J{A}^{-s} (H-\lambda \mp i \mu )^{-1} \J{A}^{-s} \phi \right\|  \\ & \quad  \leq 
c_0^{-1}
\left( 
 \left(\frac{\delta }{\tilde{R}} \right)^{1/2} + \delta (\tilde{R})^{s-1/2}\frac{ (1-s)^{1-s} (2-s)}{ (2s -1)} 
\right) ^2e^{\delta c_2\tilde{R} /c_0}
 \| \phi \| \label{1}
\end{align}
holds for all $1/2 < s \leq 1$ and $\phi \in \SCR{H}$, where $\J{\cdot} = (1+ \cdot ^2)^{1/2}$. Moreover, for all $s\geq 1$, 
\begin{align}\nn
& \sup_{| \lambda| \geq 3R, \ \mu >0} \left\| \J{A}^{-s} (H-\lambda \mp i \mu )^{-1} \J{A}^{-s} \phi \right\| \\ & \quad \leq \label{2}
c_0^{-1} \left( 
\left( \frac{\delta}{\tilde{R}} \right)^{1/2} + 2 \delta s^{-s} (s-1) ^{s-1} (\tilde{R})^{1/2}
\right) ^2e^{\delta c_2 \tilde{R} /c_0}  \| \phi\|
\end{align}
holds. In particular, if we can select $s>0$ such that it is sufficiently large and $c_2=c_2(R)$ satisfies $c_2 \tilde{R} \to 0$ as $R \to \infty$, then for some $0< \ep_{R,s} \ll 1$, 
\begin{align*}
\sup_{| \lambda| \geq 3R, \ \mu >0} \left\| \J{A}^{-s} (H-\lambda \mp i \mu )^{-1} \J{A}^{-s} \phi \right\|  \leq \ep_{R,s} \| \phi \|
\end{align*} 
holds.
\end{Thm}
Here, we define $\tilde{\varphi} \in C^{\infty} ({\bf R})$, which satisfies $\tilde{\varphi} (s) = 1 $ if $|s| \geq 5R $ and $=0$, if $|s| \leq 4R$. Then, it follows that 
\begin{align*}
& \sup_{\lambda \in {\bf R}, \, \mu >0} 
\left\| 
 \J{A}^{-s} \tilde{\varphi}(H) (H-\lambda \mp i \mu )^{-1} \tilde{\varphi} (H) \J{A}^{-s} \phi
\right\| \\ & \leq 
\sup_{ |\lambda | \leq {3 R}, \,  \mu >0} 
\left\| 
 \tilde{\varphi}(H) (H-\lambda \mp i \mu )^{-1} \tilde{\varphi} (H) \right\| \left\| \J{A}^{-s} \right\| ^2 \left\| \phi
\right\| \\ & \quad + 
\sup_{|\lambda | \geq {3 R}, \, \mu >0} 
\left\| 
 \J{A}^{-s} \tilde{\varphi}(H) \J{A}^{s'} \right\|^2 \left\| \J{A}^{-s'}  (H-\lambda \mp i \mu )^{-1} \J{A}^{-s'} \right\| \| \phi \|, 
\end{align*}
where $1/2 < s' \leq s$ with $s' \leq 1$. By noting commutator expansion (see, \S{2}), we get 
\begin{align*}
 \sup_{ |\lambda | \leq {3 R}, \,  \mu >0} 
\left\| 
 \tilde{\varphi}(H) (H-\lambda \mp i \mu )^{-1} \tilde{\varphi} (H) \right\| &\leq C R^{-1}, \\ 
 \left\| 
 \J{A}^{-s} \tilde{\varphi}(H) \J{A}^{s'} \right\|^2 &\leq 1 + CR^{-1},
\end{align*} 
which yields the LAP: 
\begin{align*}
& \sup_{ \lambda \in {\bf R},  \ \mu >0} \left\| \J{A}^{-s} \tilde{\varphi} (H) (H-\lambda \mp i \mu )^{-1} \tilde{\varphi} (H) \J{A}^{-s} \phi \right\|  \\ & \quad  \leq 
c_0^{-1} \delta
\left( 
 \left(\frac{\delta }{\tilde{R}} \right)^{1/2} + \delta (\tilde{R})^{1/2}\frac{ (1-s)^{1-s} (2-s)}{ (2s -1)} 
\right) ^2e^{\delta c_2\tilde{R}/c_0}
 \| \phi \|
\end{align*} 
for all $1/2 < s \leq 1$. As the direct consequence of Kato \cite{Ka} (see also D'Ancona \cite{D}), we obtain the following theorem:
\begin{Thm}
Under the same assumptions in Theorem \ref{T2}, 
\begin{align*}
& \frac{1}{2\pi}\int_{\bf R} \left\| 
 \J{A}^{-s} \tilde{\varphi} (H) e^{-itH} \phi 
\right\| ^2 dt \\ & \leq  c_0^{-1} \delta
\left( 
 \left(\frac{\delta }{\tilde{R}} \right)^{1/2} + \delta (\tilde{R})^{s-1/2}\frac{ (1-s)^{1-s} (2-s)}{ (2s -1)} 
\right) ^2e^{\delta c_2\tilde{R}/c_0}
 \| \phi \|^2
\end{align*}
holds for all $1/2 <s \leq 1$ and $\phi \in \SCR{H}$, where $\tilde{\varphi} \in C^{\infty} ({\bf R})$, which satisfies $\tilde{\varphi} (s) = 1 $ if $|s| \geq 5R $ and $=0$ if $|s| \leq 4R$. 
\end{Thm} 
Theorems \ref{T1} and \ref{T2} yield the following corollary:
\begin{Cor}
Let $I := \sigma (H) \cap \left( (-\infty, -5R] \cup [5R, \infty) \right) $. Then $I \subset \sigma_{\mathrm{ac}} (H)$, where $\sigma_{\mathrm{ac}}$ denotes a set of all absolutely continuous spectra of $H$.
\end{Cor}
\begin{Rem}
As usual, $c_0>0$ is not a large number. For example, consider $H= p^2$, $p= -i d/dx$, $\SCR{H}= L^2({\bf R})$ and the Yokoyama-type conjugate operator $A = (xp (1+p^2)^{-1} + (1+ p^2)^{-1} px)/2$, $x \in {\bf R}$. Formally, it then follows that 
\begin{align*}
i[H,A]^0 = 2p^2 (1+ p^2)^{-1} = 2 + \CAL{O}(R^{-1}), \quad i[[H,A]^0,A]^0 = \CAL{O}(R^{-2}).
 \end{align*}
 Hence, deducing the smallness  exclusively with \eqref{1} may be difficult. On the other hand, by taking $\tilde{R} = 4$ and $s= 4$, \eqref{2} deduces $\sup_{|\lambda| \geq 3R} \| \J{A}^{-s} (H-\lambda \mp i \mu )^{-1} \J{A}^{-s} \phi \| <1$, such that the estimate can be applied to an operator written as a form (selfadjoint operator + complex perturbation). See, e.g., Wang \cite{Wa} and the references therein.  
\end{Rem}
\begin{Rem}
Under some situations, it follows that 
$$\sup_{|\lambda| \geq 3R} \| \J{A}^{-s} (H-\lambda \mp i \mu )^{-1} \J{A}^{-s} \phi \| \to 0 , \quad \mbox{as }R \to \infty, $$ and it may be expected that the same estimate holds for the general Hamiltonian $H$ and its conjugate operator $A$. However, there are some counterexamples, e.g., by taking $H = p$ on $L^2({\bf R})$. Then, $i[H, x]^0 = 1 \geq 0$, such that for any fixed $\lambda _0$ there exists $\nu>0$, where
\begin{align} \label{4}
\lim_{\mu \to 0+} \|  \J{x}^{-s} (H - \lambda _0 -i \mu)^{-1} \J{x} ^{-s} \| = \nu
\end{align} 
holds. On the other hand, suppose that for all $\ep >0$, there exists $R=R_{\ep} \gg 1$ such that 
\begin{align} \label{3}
\lim_{\mu \to 0+} \|  \J{x}^{-s} (H - R -i \mu)^{-1} \J{x} ^{-s} \| \leq \ep
\end{align}
holds. Then, by the unitary transform $e^{i(R-\lambda _0) x}$, we notice that \eqref{3} is equivalent to 
\begin{align*}
\lim_{\mu \to 0+} \|  \J{x}^{-s} (H - \lambda _0 -i \mu)^{-1} \J{x} ^{-s} \| \leq \ep, 
\end{align*}
which contradicts \eqref{4} by taking a sufficiently small $\ep>0$ compared to $\nu$. Hence, some additional assumptions are needed with regard to $H$, $A$, and $\varphi $ in order to prove the high-energy decay property for the resolvent. (The decaying estimate with a micro local parameter has been proven by Royer \cite{Ro}.)
\end{Rem} 
With a similar approach, we can prove a smoothing-type limiting absorption principle. In this paper, we additionally assume the following statements: 
\begin{Ass}\label{A2}
Suppose Assumption \ref{A1} and let $B:= i[H,A]^0$. Then, commutators 
\begin{align*}
\mathrm{ad}_A^2 (H) := i[B, A]^0, \quad \mathrm{ad}_A^3 (H) := i[[B,A]^0,A]^0 , \quad \mathrm{ad}_A^4 (H):= 
i[\mathrm{ad}_A^3 (H),A ]^0
\end{align*}
are well-defined and can be extended to bounded operators. Moreover, for some $c_0 \leq \tilde{c} \leq c_1$, the commutator $i[H,A]^0$ can be written as 
\begin{align} \label{aad1}
i[H,A]^0 = \tilde{c} I + J + K, 
\end{align}
where $I$ is the identity operator on $\SCR{H}$, and $J$ and $K$ are bounded operators. Moreover, they satisfy the following four conditions: \\
(g). For large $R$, there exists a sufficiently small constant $\delta >0$, such that 
\begin{align*} 
\| K \varphi (H) \| \leq \delta. 
\end{align*}
(h). $J$ commutes with $H$. \\
(i). For the same $0 \leq \beta \leq 1/2$, $[\J{H}^{\beta},K ]^0 \J{H}^{\beta} A$ is a bounded operator. \\
(k). For the same $0 \leq \beta \leq 1/2$, the operator $\J{H}^{\beta} i[i[H,A]^0, A]^0 \J{H}^{\beta} $ is a bounded operator.
\end{Ass}  
\begin{Thm}\label{T3}
For the same $\beta$ in Assumption \ref{A2}, let us define 
\begin{align*}
P := \J{H}^{2\beta} A + A \J{H}^{2 \beta}. 
\end{align*}
Then, under Assumptions \ref{A1} and \ref{A2}, for all $s >1/2$, the following limiting absorption principle holds:
\begin{align} \label{Ad5}
\sup_{|\lambda| \geq 3R, \mu >0} \|  \J{P}^{-s} \J{H}^{\beta} (H - \lambda \mp i \mu )^{-1} \J{H}^{\beta} \J{P}^{-s} \phi \| \leq C \| \phi \| .
\end{align} 
Moreover, for all $1/2 < s \leq 1$, the high-energy smoothing estimates are
\begin{align*}
\int_{\bf R} \left\| 
 \J{P}^{-s} \J{H}^{\beta} \tilde{\varphi} (H) e^{-itH}\phi
\right\| ^2 dt \leq C \| \tilde{\varphi} (H) \phi \|^2, 
\end{align*}
where $\tilde{\varphi} \in C^{\infty} ({\bf R})$, which satisfies $\tilde{\varphi} (s) = 1 $, if $|s| \geq 5R $, and $=0$, if $|s| \leq 4R$. 
\end{Thm}
\section{Preliminaries}
In this section, we introduce the most important Lemma in this paper: 
\begin{Lem} \label{L1} Under Assumption \ref{A1}, the following statements hold: \\ 
(a1) $\D{A}$ reduces $\varphi(H)$, i.e., 
\begin{align*}
 \varphi (H) \D{A} \subset \D{A}. 
\end{align*}
(a2) $\| [\varphi (H) ,A] \| \leq CR^{-1} $ holds, (denoted by $[\varphi(H),A] = \CAL{O}(R^{-1})$). \\ 
(a3) Additionally, if we suppose Assumption \ref{A2}, then for all $0 \leq \beta \leq 1/2$, $ [\J{H}^{2\beta} , i[H,A]^0 ]^0 $ is a bounded operator. Moreover, $[\J{H}^{\beta}, A]^0 \J{H}^{\beta} $ is also a bounded operator.
\end{Lem}
\Proof{
We employ so-called commutator expansion:
Let $A_0$ and $B_0$ be selfadjoint operators with 
\begin{align*}
\| i[A_0,B_0]^0 \| < \infty, \quad \| \mathrm{ad} _{A_0}^{j} (B_0) \| < 0,  
\end{align*}
for some integer $2 \leq j $. For $0 \leq \rho \leq 1$, suppose $f \in C^{\infty} ({\bf R})$ satisfies $|\partial _s^k f(s)| \leq C_k \J{s}^{\rho - k}$, $k \geq 0$. Then, 
\begin{align*}
i[f(A_0), B_0] = \sum_{k=1}^{j-1} \frac{1}{k!} f^{(k)} (A_0) \mathrm{ad}_{A_0}^{k} (B_0) + R_{j} (f,A_0,B_0) 
\end{align*}
where $R_j (f,A_0,B_0)$ satisfies 
\begin{align*}
\| 
(A_0 + i)^{j-1} R_j (f,A_0,B_0)
\| \leq C (f^{(j)})\| \mathrm{ad}_{A_0}^{j} (B_0) \|.
\end{align*}
The proof of this lemma is given in Sigal--Soffer \cite{SS2} and as Lemma C.3.1 in Dere{z}i\'{n}ski and G\'{e}rard \cite{DG}. Commutator expansion with $j=2$ is such that
\begin{align*}
A \varphi(H) \J{A}^{-1} = [A,\varphi(H)] \J{A}^{-1} + \varphi(H) A \J{A}^{-1}
\end{align*}
will be a bounded operator, since $\varphi'$ and $\varphi ''$ are bounded functions. Moreover, by the construction of $\varphi$, we have $|\varphi '| , |\varphi ''| \leq \CAL{O}(R^{-1})$, and hence (a1) and (a2) are proven. Now we prove (a3): by $2 \beta \leq 1$ and Assumption \ref{A2}, $\| [B,H] ^0 \| \leq C$ with
 $B= i[H,A]^0$, we can use commutator expansion to a pair of operators, $B$ and $\J{H}^{2 \beta}$. Hence, $\| [\J{H}^{2 \beta}, B ]^0\| \leq C$. Moreover, formally we have 
\begin{align*} 
i[\J{H}^{\beta}, A ]^0 \J{H}^{\beta} = i\J{H}^{\beta} A \J{H}^{\beta} -iA \J{H}^{2 \beta} = 
\J{H}^{\beta} i[A, \J{H}^{\beta}]^0 + i[\J{H} ^{2 \beta},A] ^0. 
\end{align*}
By commutator expansion with $j=2$, we have $\J{H}^{\beta} i[A, \J{H}^{\beta}]^0 $ as a bounded operator, since $| \J{s}^{\beta} (\J{s}^{\beta} ) ' | \leq C \J{s}^{2 \beta -1}  \leq C$.
}

\section{Proof of Theorems}
We now prove Theorems \ref{T1} and \ref{T2}. The proof for Theorem \ref{T2} involves imitates the approach of \cite{Mo} (see Isozaki \cite{Is}). 
\subsection{Proof of Theorem \ref{T1}} 
Let $\lambda _0 \in (- \infty, -3R] \cup [3R, \infty)$ be an eigenvalue of $H$, and let $\psi \in \SCR{H}$ be an eigenfunction: i.e., $H \psi = \lambda _0 \psi$. Then, using Virial's theorem, Assumption (e), and $\varphi (\lambda_0) =1 $, we have 
\begin{align*}
0 = (i[\lambda_0 ,A] \varphi(H) \psi , \varphi(H)\psi ) = (i[H,A]^0 \varphi(H) \psi , \varphi(H) \psi ) \geq c_0 \| \varphi(H) \psi\| ^2= c_0 \| \psi \|^2
\end{align*}
which implies $\psi \equiv 0$.
\subsection{Proof of Theorem \ref{T2}}
For a small parameter $\ep >0$, $\mu >0$ and a large parameter $\lambda \in {\bf R}$ with $|\lambda| \geq 3R$, we define 
\begin{align*}
G(\ep) &:= (H-\lambda - i \mu - i \ep M^{\ast} M)^{-1} , \\ 
T (\ep)&:= H - \lambda -i \mu - i \ep M^{\ast} M
\end{align*}
with 
\begin{align*}
M = ( \varphi (H) i[H,A]^0 \varphi(H) )^{1/2}.
\end{align*}
The following lemma immediately holds: 
\begin{Lem} \label{L2} Under Assumption \ref{A1}, the following statements hold:\\
(b1). $T (\ep)$ is a closed operator, and for all $z \in {\bf C}_+$, $0 \in \rho(T(\ep))$, and $G(\ep)$ is analytic on $C_+$.  \\ 
(b2). Let $M_1$ satisfy $M_1^{\ast} M_1 \leq M^{\ast}M$. Then, for all bounded operators $\CAL{B}$, 
\begin{align*}
\| M_1 G(\ep) \CAL{B} \| \leq \ep^{-1/2} \| \CAL{B} G(\ep) \CAL{B} \|^{1/2}
\end{align*}
holds. \\ 
(b3). For some $0< \hat{\ep} <1$, suppose that $\ep < R^{1- \hat{\ep}}$. Then there exists $\delta =\delta (R)$, with $\delta \to 1$ as $R \to \infty$ such that $\left\| G(\ep) \right\| \leq  \delta (c_0 \ep )^{-1}$. \\ 
(b4). $\left\|(\varphi (H) -1) G(\ep)  \right\| \leq C R^{-1}$.  
\end{Lem}
\Proof{
We only need to prove (b3) and (b4). Proofs for (b1) and (b2) can be found in \cite{Mo} and \cite{CFKS} (or see the proof for Lemma \ref{L5} with $\beta =0$). We divide $G(\ep)$ into $\varphi (H) G(\ep) + (1- \varphi (H))G(\ep) $. By (b2) with $(\varphi (H)\sqrt{c_0})^2 \leq  M^{\ast}M $ and $\CAL{B}=1$, we have $\| \varphi (H) G(\ep)\| \leq (c_0 \ep)^{-1/2} \|G(\ep) \|^{1/2} $. On the other hand, by 
\begin{align*}
G(\ep) - (H-\lambda -i \mu )^{-1} = i \ep (H- \lambda -i \mu)^{-1} M^{\ast}M G(\ep)
\end{align*}
and for $|\lambda| \geq 3R$, 
\begin{align*}
\left\| 
(\varphi (H) -1 )  (H -\lambda -i \mu)^{-1}
\right\| \leq R^{-1}, 
\end{align*}
holds. Thus,
\begin{align} \label{aad17}
\left\| 
(\varphi (H) -1) G(\ep)
\right\| \leq R^{-1} + \ep c_1 R^{-1} \| G(\ep) \| . 
\end{align}
Hence, we derive inequality 
\begin{align*}
\left\| 
G(\ep)
\right\| \leq (c_0 \ep)^{-1/2 } \left\| 
G(\ep)
\right\| ^{1/2} +R^{-1} + \ep c_1 R^{-1} \left\| 
G(\ep)
\right\| . 
\end{align*} 
By assuming $\ep c_1 \leq R^{1-\hat{\ep}}$ for some $0< \hat{\ep}<1 $, we get 
\begin{align*}
\left\| 
G(\ep)
\right\| \leq \ep ^{-1}c_0^{-1} \left( 
1- \frac{\ep c_1}{R}
\right) ^{-1} \leq  \ep ^{-1}c_0^{-1} \left( 
1- R^{-\hat{\ep}}
\right) ^{-1} =  \ep ^{-1}c_0^{-1} \delta.
\end{align*}
By using (b3) for \eqref{aad17}, we immediately obtain (b4). 

}
Define 
\begin{align*}
W(\ep) := (|A| +1) ^{-s} (\ep |A| +1)^{s- 1}
\end{align*}
and 
\begin{align*}
F (\ep) := W(\ep) G(\ep) W(\ep) . 
\end{align*}
For simplicity, we denote $W(\ep) =W$, $F(\ep) =F$, $G(\ep) =G$, $i[H,A]^0 = B$ and $d/d \ep = '$. Then, 
\begin{align*}
-i F' = -i (W' GW +WGW') + WGM^2GW =: \sum_{j=1}^4 L_j 
\end{align*}
holds, where 
\begin{align*}
L_1 &= -i (W'GW + WGW'), \\ 
L_2 &= WG(\varphi -1) B (\varphi -1) GW, \\ 
L_3 &= WG (\varphi -1) BGW + WGB (\varphi -1) GW, \\ 
L_4 &= WGBGW.
\end{align*}
First, we prove that 
\begin{align}
\|  W' \| & \leq (1-s) c(s) \ep ^{s-1}, \label{5} \\ 
\| A W \| & \leq  \ep ^{s-1}, \label{8}
\end{align}
where $c(s) := (1-s)^{1-s} (2-s)^{s-2} \leq 1$, 
because 
\begin{align*}
|\ |A| (|A| +1) ^{-s} (\ep |A| +1) ^{s -2}  \| &= \sup_{\lambda \geq 0} \lambda (\lambda +1)^{-s} (\ep \lambda +1 )^{s -2}  \\ &= 
\ep^{s -1} \sup_{t \geq 0} t^{1-s} (t+ 1)^{s-2}
\end{align*}
provides \eqref{5}. Second, we estimate $L_j$, $j=1,2,3,4$. Here, we note that 
$$c_0 \varphi(H) ^2 = (c_0 ^{1/2} \varphi(H) )^{\ast} (c_0^{1/2} \varphi(H) ) \leq \varphi(H) i[H,A]^0 \varphi(H) 
$$ and Lemma \ref{L2}(b2) provides 
\begin{align*}
\|\varphi(H) GW \| \leq c_0^{-1/2} \ep^{-1/2} \|WGW \| ^{1/2}= (c_0 \ep)^{-1/2} \| F \|^{1/2}. 
\end{align*}
By Lemma \ref{L2}(b4),
\begin{align}
\label{7} \| (1-\varphi (H)) GW  \| \leq CR^{-1}
\end{align}
holds, and thus we have 
\begin{align}
\| GW \| \leq (c_0 \ep)^{-1/2} \| F \|^{1/2} + \CAL{O}(R^{-1}).  \label{6}
\end{align}
By \eqref{5} and \eqref{6}, we have
\begin{align*}
\| L_1 \| \leq 2(1-s) c(s) \ep ^{s-1} (\CAL{O}(R^{-1}) + \|F /(\ep c_0)  \|^{1/2}). 
\end{align*}
By \eqref{7} and $\| B \| \leq C $, we have
\begin{align*}
\| L_2 \| \leq \CAL{O}(R^{-2}). 
\end{align*}
By \eqref{7}, \eqref{6}, and $\|B \| \leq C$, we have  
\begin{align*}
\| L_3 \| \leq \CAL{O}(R^{-1}) (\CAL{O}(R^{-1}) + \|F / ( \ep c_0) \|^{1/2}) = \CAL{O}(R^{-2}) + \CAL{O}(R^{-1}) \| F/ (\ep c_0 ) \|^{1/2}
\end{align*} 
We divide $-i L_4$ into $L_5 + L_6$ with 
\begin{align*}
L_5 := WAGW -WGAW, \quad L_6 := i \ep WG [M^2, A] GW . 
\end{align*}
By \eqref{8} and \eqref{6}, 
\begin{align*}
\| L_5 \| \leq 2c(s) \ep^{s-1} (\CAL{O}(R)^{-1} + \| F/ ( \ep c_0 ) \|^{1/2}).
\end{align*}
By Lemma \ref{L1}(a2) and 
\begin{align*}
[M^2,A ] = [\varphi(H), A] i[H,A] \varphi(H) + \varphi(H) i[[H,A],A] \varphi(H) + \varphi(H)i[H,A] [\varphi(H),A ], 
\end{align*}
we have 
\begin{align*}
\| L_6 \| & \leq \ep \left( 
\CAL{O}(R^{-1}) + \| F/ ( \ep c_0)  \|^{1/2}
\right) ^2 \left( 
\CAL{O}(R^{-1} ) + c_2 
\right)  \\ & \leq 
((c_2/c_0)  + \CAL{O}(R^{-1}) ) \| F \| + \ep^{1/2} \CAL{O}(R^{-1}) \|  F \|^{1/2} + \ep \CAL{O}(R^{-2}).
\end{align*}
Then, for some $\delta = \delta _R >1$ with $\delta \to 1$ as $R \to \infty$, 
\begin{align*}
\| F' \| & \leq  \sum_{j=1}^6 \| L_j \|   \leq 2\delta (2-s) c(s) c_0^{-1/2} \ep ^{s-3/2} \| F \|^{1/2} + (\delta c_2 /c_0 ) \| F \| + \CAL{O}(R^{-1})
\end{align*}
holds. By following the argument in \cite{Mo}, for all $0< \tilde{R} \leq R^{1-\hat{\ep}}$ and for some $\delta$  it holds that
\begin{align*}
F(\ep) & \leq \left(  ( \delta F(\tilde{R}) )^{1/2}
+ \frac12 \int_{\ep}^{\tilde{R}} \left( 2 \delta (2-s) c(s)c_0^{-1/2} t^{s-3/2} \right) e^{-\delta  c_2 (\tilde{R}-t)/(2 c_0)} dt 
\right) ^2 e^{\delta c_2 \tilde{R} /c_0} \\ & \leq c_0^{-1}
\left( 
 \left( \frac{\delta }{\tilde{R}} \right)^{1/2} + \frac{2\delta (2-s)(\tilde{R}^{s-1/2})}{2s-1} c(s) 
\right) ^2e^{\delta \tilde{R} c_2/c_0},
\end{align*}
which proves \eqref{1}. Here, we remark that 
\begin{align*}
(2-s) c(s) (2s -1)^{-1} = (1-s)^{1-s} (2-s) (2s -1)^{-1} \geq 1.
\end{align*}

Next we prove \eqref{2}. Except for replacing $W(\ep) $ to $W := (1+|A| )^{-s} $, the approach is the same as in the proof for \eqref{1}. The important estimate is 
\begin{align*}
\| W' \|  = 0 , \quad \mathrm{and} \quad 
\| W A \| \leq s^{-s} (s-1) ^{s-1} = \left( 
1- \frac{1}{s}
\right) ^{s} (s-1)^{-1}=: C(s). 
\end{align*}
Then, it follows that 
\begin{align*}
\|  L_1\| \leq 0, \quad \mbox{and} \quad \| L_5 \| \leq 2C(s) (\CAL{O}(R)^{-1} + \| F/ ( \ep c_0 ) \|^{1/2})
\end{align*}
Hence, we have 
\begin{align*}
\| F' \| \leq 2 \delta C(s) c_0^{-1/2} \ep ^{-1/2} \| F \|^{1/2} + (\delta c_2 /c_0) \| F \| + \CAL{O}(R^{-1}). 
\end{align*}
Let $\tilde{R}$ be a positive constant with $\tilde{R} < R$ and then 
\begin{align*}
F(\ep) &\leq \left( 
(\delta F(\tilde{R}) ) ^{1/2} + \frac{1}{2} \int_{\ep}^{\tilde{R}} \left(  2 \delta C(s) c_0^{-1/2} t^{-1/2} \right) e^{- \delta c_2(\tilde{R} -t) /(2c_0) } dt 
\right) ^2 e^{\delta c_2 (\tilde{R} - \ep) /c_0} \\ & \leq 
c_0^{-1} \left( 
\left( \frac{\delta}{\tilde{R}} \right)^{1/2} + 2 \delta C(s) (\tilde{R})^{1/2}
\right) e^{\delta c_2 \tilde{R} /c_0} 
\end{align*}
holds, which gives \eqref{2}. 

\section{Smoothing estimate with high-energy cut-off}
In this section, we prove Theorem \ref{T3}. In the case where $\varphi \in C_0^{\infty} ({\bf R})$, a more generalized estimate was obtained by M\o ller-Skibsted \cite{MS}. For $0 \leq \beta \leq 1/2$, we define 
\begin{align*}
G_{\beta}(\ep) &:= (H - \lambda -i \mu - i \ep \varphi (H) \J{H}^{\beta} i[H,A]^0 \J{H}^{\beta} \varphi (H)  ) ^{-1} \\ 
T_{\beta}(\ep) &:= H - \lambda -i \mu - i \ep \varphi (H) \J{H}^{\beta} i[H,A]^0 \J{H}^{\beta} \varphi (H) , \\ 
M_{\beta} &:= ( \varphi (H) \J{H}^{\beta} i[H,A]^0 \J{H}^{\beta} \varphi (H) )^{1/2}, 
\end{align*} 
where we always assume that $\mu >0 $ and $\ep >0$ are sufficiently small. For simplicity, we may denote $$\J{H}^{\beta} \varphi(H) = g_{\beta}(H) =g(H) , \quad M_{\beta} = M. $$
\begin{Lem}\label{L4}
For any fixed $\ep$ and $\mu$, we define $H(\ep):=\J{H}^{\beta} G(\ep) \J{H}^{\beta}$. Then, $H(\ep)$ satisfies 
\begin{align*}
\left\| 
H(\ep) 
\right\| \leq C (\mu^2/\lambda ^{4 \beta} + \ep ^2) ^{-1/2}. 
\end{align*}
\end{Lem}
\Proof{
Here, let us consider the operator 
\begin{align*}
\hat{T} := \J{H}^{-\beta} (H - \lambda -i\mu -i \ep M^{\ast} M ) \J{H}^{-\beta}. 
\end{align*}
Clearly, $\hat{T}$ can be defined on $\D{\J{H}^{1-2 \beta}}$. Then, for all $u\in \D{\J{H}^{1-2 \beta}} $ with $\| u \| = 1$,
\begin{align*}
\left\| 
\hat{T} u
\right\|^2 &= \left\| \J{H}^{-\beta}(H-\lambda )\J{H}^{-\beta} u \right\|^2 + \left\| \J{H}^{-\beta}
(\mu + \ep M^{\ast} M) \J{H}^{-\beta} u 
\right\| ^2 \\ & \quad  - \left( 
i\ep [\J{H}^{-\beta} (H-\lambda ) \J{H}^{-\beta}, \J{H}^{-\beta} M^{\ast} M \J{H}^{-\beta}   ] u, u
\right) \\ 
 & \geq I_1 + I_2 + I_3,  
\end{align*}
where 
\begin{align*}
I_1 &:=  \left\| \J{H}^{-\beta}(H-\lambda )\J{H}^{-\beta} u \right\|^2 , \\ 
I_2 &:=   \mu^2\| \J{H}^{-2 \beta} u \|^2 + \mu \ep c_0 \| \varphi (H)u \|^2 + \ep^2 c_0^2  \| \varphi (H)^2u \|^2 , \\
I_3 &:= -\ep \left( 
i[\J{H}^{-\beta} (H-\lambda ) \J{H}^{-\beta}, \varphi(H) i[H,A]^0 \varphi (H)]  u,u
\right)   . 
\end{align*}
By \eqref{aad1} and the boundedness of $\left\| \J{H}^{-\beta} (H-\lambda ) \J{H}^{-\beta} u \right\| $, we get
\begin{align*}
I_3 &\leq - \ep    \left( 
i[\J{H}^{-\beta} (H-\lambda ) \J{H}^{-\beta}, \tilde{c}\varphi(H) ^2 +J\varphi(H)^2 ]  u,u
\right) \\ & \quad + 2 \ep \mathrm{Im} \left( 
{K} \varphi (H) u, \J{H}^{-\beta} (H-\lambda ) \J{H}^{-\beta}  \varphi(H) u
\right) . 
\end{align*}
By the assumption of $J$ and $K$, we have 
$$\| {K} \varphi (H) \| \leq \delta, \quad \mbox{and} \quad i[\J{H}^{-\beta} (H-\lambda ) \J{H}^{-\beta},  \tilde{c}\varphi(H) ^2 +J \varphi(H)^2 ] =0, $$ 
and we also have 
\begin{align*}
\left| I_3 \right| \leq & 4 \ep \delta \| \J{H}^{-\beta} (H-\lambda ) \J{H}^{-\beta} (1- \varphi (H)  ) \varphi (H) u \| \\ & \quad + 4 \ep \delta \| \J{H}^{-\beta} (H-\lambda ) \J{H}^{-\beta}\varphi (H)^2 u \| . 
\end{align*}
Here, we define
\begin{align*}
I_4 &:= \left\| \J{H}^{-\beta}(H-\lambda )\J{H}^{-\beta} (1-\varphi (H)^2 ) u  \right\| ^2 
\\ & \geq  \left\| \J{H}^{-\beta}(H-\lambda )\J{H}^{-\beta} (1-\varphi (H) ) u  \right\| ^2 \\ 
I_5 &:= \left\| \J{H}^{-\beta}(H-\lambda )\J{H}^{-\beta} \chi (|H| \leq 2 \lambda ) \varphi (H)^2 u  \right\| ^2, \\ 
I_6 &:=\left\| \J{H}^{-\beta}(H-\lambda )\J{H}^{-\beta} \chi (|H | \geq 2\lambda ) u  \right\| ^2, \\ 
I_7  &:=  \left\| \J{H}^{-\beta} (H- \lambda) \J{H}^{-\beta} {\varphi (H)(1-\varphi (H)) } u  \right\|^2,  
\end{align*}
where we use $\varphi (H)  \chi (|H | \geq 2\lambda )= \chi (|H | \geq 2\lambda ) $. Then, we have $I_1 \geq I_4 + I_5 + I_6   $, $I_3 \leq 4 \ep \delta (\sqrt{I_5} + \sqrt{I_6} + \sqrt{I_7}) $ and 
\begin{align*}
& I_5/2 - 4 \ep \delta \sqrt{I_5} \geq - 16 \ep ^2 \delta ^2 =  - 16 \ep ^2 \delta ^2 \| u \|^2 \\ & \geq 
-64 \ep^2 \delta ^2 \left( 
\| (1- \varphi (H) ) u \|^2 + \left\| \varphi (H) \chi (|H| \leq 2 \lambda) u \right\|^2 + \left\|  \chi (|H| \geq 2 \lambda) u \right\|^2
\right). 
\end{align*}
Hence, we have 
\begin{align*}
I_1-|I_3| & \geq I_4 + I_5/2 +I_6 -4 \ep \delta \left( 
\sqrt{I_6} + \sqrt{I_7}
\right) \\ & \quad -64 \ep^2 \delta ^2 \left( 
\left\| 
(1-\varphi (H)) u
\right\|^2 + \left\| 
\varphi (H) \chi (|H| \leq 2 \lambda) u
\right\|^2 + \left\| 
\chi (|H| \geq 2 \lambda) u
\right\|^2
\right). 
\end{align*}
We thus divide into 
\begin{align*}
I_1+I_2+I_3  \geq J_1+J_2 +J_3, 
\end{align*}
where 
\begin{align*}
J_1 &:= I_6-4 \ep \delta \sqrt{I_6} -64 \ep^2 \delta ^2 \left\| \chi (|H| \geq 2 \lambda) u \right\|^2, \\ 
J_2 &:= I_4/2 +I_5/2 -4 \ep \delta \sqrt{I_7} -64 \ep^2 \delta ^2 \left\| (1- \varphi (H)) u \right\|^2, \\ 
J_3 &:= I_4/2 + I_2 -64 \ep^2 \delta ^2 \left\| \varphi (H)\chi (|H| \leq 2 \lambda) u \right\|^2
\end{align*}
First, we estimate $J_1$. It is clear that for large $|y|\geq 2 \lambda $
\begin{align*}
 \left| (y-\lambda) \J{y}^{-2 \beta} \right| \geq \frac12,
\end{align*}
and hence by taking $R>0$ to be sufficiently large, we have 
\begin{align} \label{aad6}
 & \left\| \J{H}^{-\beta}(H-\lambda )\J{H}^{-\beta} \chi (|H | \geq 2\lambda ) u  \right\| \geq \left\| \chi (|H | \geq 2\lambda ) u  \right\|/4 \end{align}
which gives 
\begin{align} \nn
& I_6 - 64 \ep ^2 \delta ^2 \| \chi(|H| \geq 2 \lambda) u \|^2 - 4 \ep \delta \sqrt{I_6} \geq (1- 1028 \ep ^2 \delta ^2 ) I_6 - 4 \ep \delta \sqrt{I_6} \\ 
 &= v_1 \left( 
 \sqrt{I_6} - \frac{2\ep\delta}{ v_1} \right) ^2 - \frac{4\delta ^2\ep^2}{ v_1},  \label{aad5}
\end{align}
where we assume $v_1 := (1-1028 \ep^2 \delta ^2) \geq 1/2$. If $\chi (|H| \geq 2 \lambda) u\equiv 0$, then $\eqref{aad5} = 0 = \| \chi (|H| \geq 2 \lambda) u \|^2 $, and if $\chi (|H| \geq 2 \lambda) u \neq 0$, \eqref{aad6} implies $ \sqrt{I_6} \geq 1/4$ then for some certain $\ep$-independent positive constant $v_2 >0$, 
\begin{align*}
\eqref{aad5} \geq v_2 = v_2 \| \chi (|H| \geq 2 \lambda) u \|^2
\end{align*}
holds, which gives 
\begin{align*}
& J_1 \geq v_2 \| \chi (|H| \geq 2 \lambda) u \|^2.
\end{align*}
Next we estimate $J_2$. By the definition of $\varphi$, $\sqrt{I_7}$ can be divided into 
\begin{align*}
\sqrt{I_7} &= \left\|  \J{H}^{-\beta}(H-\lambda )\J{H}^{-\beta} \varphi(H) (1-\varphi(H))\chi (3R/2 \leq |H | \leq 2\lambda ) u \right\| \\ & \quad + 
\left\|  \J{H}^{-\beta}(H-\lambda )\J{H}^{-\beta} \varphi(H) (1-\varphi(H))\chi (|H| \leq 3R/2 ) u \right\| \\ & =: I_8+I_9.
\end{align*}
On the support of $(1-\varphi (H) )$, 
\begin{align*} 
& \left\| \J{H}^{-\beta} (H -\lambda) \J{H}^{- \beta}(1-\varphi (H) )u \right\| \geq R(1+4R^2)^{- \beta}
\left\| 
(1-\varphi (H) )u
\right\|  \\ &\geq 5^{-\beta} R^{1-2 \beta}  \left\| 
(1-\varphi (H) )u
\right\|
\end{align*} 
holds. Hence, there is a positive constant $v_3 \geq 1/4$ such that 
\begin{align}
& \label{aad7} I_4/2 -64 \ep^2 \delta ^2 \left\| (1- \varphi (H)) u \right\|^2 - 4 \ep \delta I_9
\geq v_3I_4- 4 \ep \delta I_9 \\ \nn & \geq 
v_3 \left\| 
 \J{H}^{-\beta}(H-\lambda )\J{H}^{-\beta}(1-\varphi(H))\chi (|H| \leq 3R/2 ) u 
\right\|^2 \\ \nn  & \quad - 4 \ep \delta  \left\| 
 \J{H}^{-\beta}(H-\lambda )\J{H}^{-\beta}(1-\varphi(H))\chi (|H| \leq 3R/2 ) u 
\right\|
\end{align}
holds. On the support of $\chi (|H| \leq 3R/2)$, $(1-\varphi (H))$ is not always $0$, since $\varphi '(s) >0$ on $R<|s| <2R$, which gives $ \eqref{aad7} \geq 0$, by selecting $\ep \delta >0$ enough small. In the same manner, we get $I_5/2-4\ep \delta I_8 \geq 0$. Thus, we have $J_2 \geq 0$. These yield 
\begin{align*}
I_1+I_2+I_3 \geq v_2  \| \chi (|H| \geq 2 \lambda) u \|^2 + J_3 
\end{align*}
Since $J_3$ is greater than 
\begin{align*}
&\frac14 \left\| \J{H}^{-\beta}(H-\lambda )\J{H}^{-\beta} (1-\varphi (H)^2 ) u  \right\| ^2 + \mu^2 \left\| \J{H}^{-2 \beta} \varphi(H) \chi(|H| \leq 2 \lambda) u \right\|^2 \\ & \quad + \left( \ep^2 c_0^2 \left\|  \chi (3R/2 \leq |H| \leq 2 \lambda) \varphi (H)^2 u \right\| ^2 -64 \ep^2 \delta ^2 \left\| 
 \chi (3R/2 \leq |H| \leq 2 \lambda) \varphi (H) u
\right\| ^2 \right) \\ & \qquad + \Big( \frac14 \left\| \J{H}^{-\beta}(H-\lambda )\J{H}^{-\beta} (1-\varphi (H)^2 )\chi (|H| \leq 3R/2) u  \right\| ^2 \\ & \quad \qquad  -64 \ep^2 \delta ^2 \left\| 
 \chi (|H| \leq 3R/2) \varphi (H) u
\right\| ^2  \Big)  \\ & \geq 
v_4  \left\| (1-\varphi (H)^2 ) u  \right\| ^2 + \mu^2 \J{\lambda} ^{-4\beta} \left\| 
\varphi (H) ^2\chi(|H| \leq 2 \lambda) u
\right\|^2 \\ & \quad + v_5 \ep^2 c_0^2 \left\|  \chi (3R/2 \leq |H| \leq 2 \lambda) \varphi (H)^2 u \right\| ^2 + 
v_6   \left\|  \chi (|H| \leq 3R/2) u \right\| ^2 \\ & \geq 
v_4 \left\| (1-\varphi (H)^2 ) u  \right\| ^2 + v_7(\ep ^2 + (\mu/|\lambda|^{2\beta} )^2   )
\left\| 
\varphi (H) ^2\chi(|H| \leq 2 \lambda) u
\right\|^2
\end{align*}
for some $(\ep,\mu,\lambda)-$independent positive constants $v_4$, $v_5$, $v_6$, and $v_7$. Finally, we obtain 
\begin{align}
I_1+I_2+I_3 \nn & \geq 
v_4 \left\| (1-\varphi (H)^2 ) u  \right\| ^2 + v_7(\ep ^2 + (\mu/|\lambda|^{2\beta} )^2   )
\left\| 
\varphi (H) ^2\chi(|H| \leq 2 \lambda) u
\right\|^2 \\ & \quad +v_2 \| \chi (|H| \geq 2 \lambda) u \|^2 \nn 
\\ & \geq C(\ep ^2 + (\mu/|\lambda|^{2\beta} )^2 \left\| u \right\|^2
 \label{aad2}
\end{align}
where we take $0<\delta ^2 \ll c_0^2$ and we use $\varphi ^2 \leq \varphi$. Inequality \eqref{aad2} entails that there exists a strict inverse operator of $\hat{T}$:
\begin{align*}
\hat{T}^{-1} = \J{H}^{\beta} (H-\lambda -i \mu -i \ep M^{\ast} M)^{-1} \J{H}^{\beta}
\end{align*}
satisfying 
\begin{align*}
\| \hat{T} \| \leq C ((\mu/|\lambda|^{2\beta})^2 + \ep ^2 )^{-1/2}.
\end{align*}
}
\begin{Lem}\label{L5}
The following statements hold: \\ 
(c1). Operator $H(\ep) :=\J{H}^{\beta} G_{\beta}(\ep) \J{H}^{\beta}$ is differentiable in $\ep$ and satisfies $$H'(\ep) = i \J{H}^{\beta} G_{\beta}(\ep) M^{2}_{\beta} G_{\beta}(\ep) \J{H}^{\beta} .$$ 
(c2). Let $Q_1$ be an operator satisfying $\D{\J{H}^{\beta}} \subset \D{Q_1 }$ and $Q_1^{\ast} Q_1 \leq M_{\beta}^{\ast} M_{\beta}$. Then, for all bounded selfadjoint operators $\CAL{B}$, 
\begin{align} \label{Ad6}
\| Q_1 G_{\beta}(\ep) \J{H}^{\beta} \CAL{B} \| \leq \ep ^{-1/2} \| \CAL{B} \J{H}^{\beta}G_{\beta}(\ep) \J{H}^{\beta} \CAL{B} \|^{1/2}
\end{align}
holds. \\
(c3). $\| G_{\beta} (\ep) \| \leq  \| \J{H}^{\beta} G_{\beta} (\ep) \J{H}^{\beta}  \| \leq C \ep ^{-1}$. \\ 
(c4). $\| \J{H}^{\beta} (1- \varphi (H)) G_{\beta}(\ep) \J{H}^{\beta}  \| \leq C$.
\end{Lem}
\Proof{
The resolvent formula 
\begin{align*}
\J{H}^{\beta} G_{\beta}(\ep) \J{H}^{\beta} -\J{H}^{\beta} G_{\beta}(\ep') \J{H}^{\beta} = i(\ep - \ep ')\J{H}^{\beta} G_{\beta}(\ep) M^2_{\beta} G_{\beta}(\ep ') \J{H}^{\beta}
\end{align*}
gives (c1). By a similar calculation and $\mu >0$, we have 
\begin{align*}
(i/2\ep)\CAL{B} \J{H}^{\beta} \left( 
G_{\beta}(\ep)^{\ast} - G_{\beta}(\ep)
\right)   \J{H}^{\beta} \CAL{B} &\geq \CAL{B} \J{H}^{\beta}
G_{\beta}(\ep)^{\ast} M^{\ast}_{\beta} M_{\beta}  G_{\beta}(\ep)
  \J{H}^{\beta} \CAL{B} \\ & \geq 
   \CAL{B} \J{H}^{\beta}
G_{\beta}(\ep)^{\ast} Q^{\ast}_1 Q_1  G_{\beta}(\ep)
  \J{H}^{\beta} \CAL{B},  
\end{align*}
which proves (c2), where we remark that 
\begin{align*}
\|  Q_1  G_{\beta} (\ep) \J{H}^{\beta} \CAL{B} \| \leq 
\| Q_1 \J{H}^{- \beta} \| \| \CAL{B} \| \| H(\ep) \|  
\end{align*}
is well defined. By \eqref{Ad6}, 
\begin{align*}
\| g(H) G_{\beta} (\ep) \J{H} ^{\beta}  \| \leq C \ep ^{-1/2} \| \J{H}^{\beta}  G_{\beta} (\ep) \J{H}^{\beta}  \|^{1/2}
\end{align*}
holds. Since $(1-\varphi ) \in C_0^{\infty}({\bf R})$, $|(1-\varphi (a))(a -\lambda -i \mu)^{-1}\J{a}^{2 \beta} | \leq C $ and 
\begin{align}
G_{\beta}(\ep) -(H- \lambda -i \mu)^{-1} = (H-\lambda-i \mu)^{-1} (i \ep M_{\beta} ^2) G_{\beta} (\ep) \label{Ad7}
\end{align} we also have 
\begin{align*}
\| (1- \varphi (H)) \J{H}^{\beta} G_{\beta} (\ep) \J{H} ^{\beta} \| \leq C (1 + \ep \| H(\ep)\| ).
\end{align*} 
Consequently we get $\| H(\ep)\| \leq C(1 + \ep ^{-1/2} \| H(\ep) \|^{1/2} + \ep \| H(\ep) \|  ) $, which proves (c3) and (c4). }

\subsection{Proof of Theorem \ref{T1}}
Now, we prove \eqref{Ad5}. The fundamental approach is exactly the same as in the proof of \eqref{1}. We define 
\begin{align*}
P &=P_{\beta } := \J{H}^{2 \beta } A +A\J{H}^{2\beta  } ,\\
W &=W_{\beta}(\ep) := (1+ | P | )^{-s} (1+ \ep |P| )^{s-1}, \quad  \\ 
G &={G}_{ \beta} (\ep) := (H- \lambda -i\mu -i \ep\varphi (H)\J{H}^{\beta} i[H,A]^0 \varphi (H) \J{H}^{\beta} )^{-1}, \\ 
F &= F_{\beta} (\ep) := W_{\beta} (\ep) \J{H}^{\beta} G_{\beta} (\ep) \J{H}^{\beta} W_{\beta} (\ep), \\ 
T & = T_{\beta} (\ep) := H - \lambda -i \mu -i \ep \varphi (H) \J{H}^{\beta} i[H,A]^0 \J{H}^{\beta} \varphi (H)
\end{align*} 
For simplicity, we use notations $F_{\beta}(\ep) =F$, $G_{\beta}(\ep )=G$, $T_{\beta} (\ep) =T$, $\J{H}^{\beta} =Z$, $B= i[H,A] ^0$, $g= g(H) = Z \varphi (H)$, $M = (\varphi(H) ZBZ \varphi (H) )^{1/2}$ and $d/d \ep ='$. A straightforward calculation shows
\begin{align*}
-i F' = -i ( W'ZGZW + WZGZW') + WZGM^2GZW =: \sum_{j=1}^4 L_j, 
\end{align*}
where 
\begin{align*}
L_1 &=  -i ( W'ZGZW + WZGZW'), \\ 
L_2 &= WZG(\varphi -1) ZBZ(\varphi -1) GZW, \\ 
L_3 &= WZG(\varphi -1) ZBZ GZW + WZGZBZ(\varphi -1)  GZW, \\ 
L_4 &= WZGZBZGZW =: L_5 +L_6 \\ 
L_5 &= i ( WZGZTAZGZW  - WZGZAT ZG ZW), \\ 
L_6 &= i\ep WZG Zi[gBg,A]Z GZW.
\end{align*}
By \eqref{5}, Lemma \ref{L5}(c4) and  
\begin{align*}
\| \varphi (H) ZGZW  \| \leq C \ep^{-1/2} \|F \|^{1/2}  
\end{align*}
we have 
\begin{align*}
\| L_1\| \leq C \ep^{s-1} ( 1 + \ep^{-1/2} \| F \|^{1/2} ).
\end{align*}
By Lemma \ref{L5}(c4), 
\begin{align*}
\| L_2 \| \leq C
\end{align*}
holds. By Lemma \ref{L5}(c4), we have  
\begin{align*}
\| WZG(1-\varphi) Z B gGZW \| \leq C \| \varphi(H) ZGZW  \| , 
\end{align*}
which yields
\begin{align*}
\|  L_3 \| \leq C (1 + \ep^{-1/2} \|F \|^{1/2}).
\end{align*}
Now we estimate $\| L_4 \|$. By \eqref{aad1}, we have 
\begin{align*}
T &= H -\lambda -i \mu -i\ep gBg \\ 
&= 
H -\lambda -i \mu - i\ep \tilde{c}g^2 -i \ep Jg^2 -i \ep gKg  
\end{align*} 
Hence, 
\begin{align*}
Z T = TZ -i \ep g[Z,K]g
\end{align*}
holds and this yields 
\begin{align*}
WZGZT AZGZW &= WZ^2A\cdot ZGZW -i \ep WZGZ \cdot Z^{-1}g[Z,K]gA\cdot ZGZW .
\end{align*}
By Assumption \ref{A2} condition (i), we have 
\begin{align*} 
\left\| Z^{-1} g[Z,K]gA \right\| \leq \|[Z,K]ZA \| +C \leq C .
\end{align*}
By $2 Z^2 A = AZ^2 + Z^2A +[Z^2,A] = P + \mbox{(bounded operator)}$, $WZ^2A = WP/2 + \mbox{(bounded operator)} $ is bounded. Consequently, we get 
\begin{align*}
\| L_5 \| &\leq C\ep (1 + \ep ^{-1/2} \| F \|^{1/2} )^2  + C s^{s-1}  (1 + \ep ^{-1/2} \| F \|^{1/2} ) \\ & \leq 
C (1+ \| F \| + s^{s-3/2} \| F \|^{1/2} ).
\end{align*}
Lemma \ref{L1} (a3) yields
\begin{align*}
\| \J{H}^{\beta}  [\J{H}^{\beta} ,A ] \| \leq C. 
\end{align*}
Moreover, by Assumption \ref{A2} condition (k),  
\begin{align*}
\| g[B,A]g \| \leq C 
\end{align*}
holds. Hence, we get 
\begin{align*}
\| L_6 \| \leq C(1 + \| F\|). 
\end{align*}
Consequently, we get \eqref{Ad5}.

\section{Application to Schr\"{o}dinger-type operators}
We now apply the main theorem to dissipative operators. Unfortunately, complex and long calculations are needed in order to check all of the conditions stated in Assumption \ref{A2}. Thus, we here only provide a sketch of the calculation. In what follows, let $| \lambda | \gg 1$, $\SCR{H} = L^2({\bf R}^n)$,  $p= -i \nabla $ on ${\bf R}^n$, and where $V$ is smooth and satisfies the condition that for all multi-index $\alpha $ and for some $\rho >0$, there exists a constant $C_{\alpha } >0$ such that 
\begin{align*}
\left| \J{x}^{\rho + |\alpha|} \partial ^{\alpha}_x V(x) \right| \leq C_{\alpha}
\end{align*}
Let us consider the Schr\"{o}dinger-type operators 
\begin{align} \label{aad12}
H = a (p^2 + b ^2)^{\gamma} +V, 
\end{align}
where $a>0$, $b \in {\bf R}$, and $ \gamma \geq 1/2$, and its conjugate operator 
\begin{align} \label{aad13}
 A = (1+ (p^2 + b ^2))^{- \gamma} p \cdot x + x \cdot p (1+ (p^2 + b ^2))^{- \gamma}.
\end{align} By simple calculation, formally, we have 
\begin{align*}
i[H,A] &= 4 \gamma a p^2 (p^2+ b ^2)^{\gamma-1} (1+ (p^2+b ^2))^{-\gamma}  \\ & \quad  -x \cdot \nabla V \times (1+p^2 +b^2)^{-\gamma} + (\mbox{ similar terms }),
\end{align*}
\begin{align*}
i[i[H,A],A] &= 16 \gamma a p^2 (p^2+b^2)^{\gamma -1} (1+p^2 +b^2)^{-2 \gamma} + (\mbox{ similar terms }) \\ & \quad + x^2 (\Delta V) (1+p^2+b^2)^{-2 \gamma} + (\mbox{ similar terms })
\end{align*}
and 
\begin{align*}
i[H,i[H,A]] &= -8 \gamma a \nabla V \cdot p (p^2+ b ^2)^{-\gamma-1} (1+ (p^2+b ^2))^{-\gamma} + ( \mbox{ similar terms }) \\ & \quad + 2a \gamma p \cdot \nabla V  \times (1+p^2 +b^2)^{-\gamma} + (\mbox{ similar terms }). 
\end{align*} 
Here, let $R$ be sufficiently larger than $C_{\alpha}$ with $a$, $b$, and $\gamma$ as constants. Then, on the support of $\varphi (H)$, it follows that $ H_0  = H-V $ is sufficiently large. Now we prove \eqref{aad1} and the four statements in Assumption \ref{A2}. 
A straightforward calculation gives 
\begin{align*}
& p^2 (p^2+ b ^2)^{-\gamma-1} (1+ (p^2+b ^2))^{-\gamma} \\ &= (p^2+b^2)^{\gamma} (1+ p^2+ b^2)^{- \gamma} -b^2 (p^2 +b^2) ^{\gamma -1} (1+p^2+b^2)^{-\gamma}  \\ 
& = 1+ \CAL{O}(( (p^2+b^2)^{\gamma})^{-1/\gamma} ) \\ 
&= 1 + \CAL{O} ({H^{-1/\gamma}}) + \CAL{O}(V (p^2 +b^2)^{-1 - 1/\gamma } ),  
\end{align*}
and hence we can divide $i[H,A]$ into $\tilde{c}I + J + K$ with $\tilde{c}= 4 \gamma a $, $J=J(H) = \CAL{O}(H^{-1/\gamma})$ and $K = \CAL{O}(V(p^2 +b^2)^{-1 - 1/\gamma } )-x\nabla V (1+p^2+b^2)^{- \gamma} + (\mbox{ similar terms })$. Since $V$ is smooth, $i[H, K]$ can be written as $\nabla V \cdot p \times \CAL{O} (1+p^2+b^2)^{-1}$. Hence, (g), (h), and (i) can be proven. Moreover, by $i[i[H,A],A] = \CAL{O}((1+p^2+b^2)^{-\gamma}) $, we also have (k). Hence, we have the following smoothing estimate: 
\begin{align} \label{aad15}
\int_{\bf R} \left\| 
\J{P}^{-s} \J{H}^{\beta} \tilde{\varphi} (H) e^{-itH} \phi 
\right\|^2 dt \leq C \| \phi \| ^2.  
\end{align}
Let $\theta \in {\bf R}$ be a constant given later. Inequality  
\begin{align*}
& \| \J{x}^{-s} \J{p} ^{\theta} \tilde{\varphi }(H)e^{-itH} \phi \| \\ & \leq 
\left\| 
\J{x}^{-s} \J{p} ^{\theta} \J{H}^{-\beta}\J{P}^s 
\right\| \| \J{P}^{-s} \J{H} ^{\beta} \tilde{\varphi} (H)e^{-itH} \phi \|
\end{align*}
gives the following smoothing estimates: 
\begin{Cor}
Suppose $H$ and $A$ are written in the form in \eqref{aad12} and \eqref{aad13}, respectively, and $\gamma \geq 1/2$. Let
\begin{align*} 
\theta = 
\begin{cases}
-4 \beta \gamma s + 2 \beta \gamma -s + 2 \gamma s &  \beta \geq (2 \gamma -1)/4\gamma, \\ 
2\beta \gamma &  \beta \geq  (2 \gamma -1)/4\gamma
\end{cases} 
\end{align*}
Then, 
\begin{align} \label{aad16}
\int_{\bf R} \left\| 
\J{x}^{-s} \J{p} ^{\theta} \tilde{\varphi} (H) e^{-itH} \phi 
\right\|^2 dt \leq C \| \phi \|^2
\end{align}
holds.
\end{Cor}
\Proof{
Let us define 
$$
\beta ^{\ast} := \frac{2 \gamma -1}{4 \gamma}. 
$$ and prove for $\beta \geq \beta ^{\ast}$ 
\begin{align} \label{aad9}
\left\| 
P \J{p}^{-4 \beta \gamma -1 + 2 \gamma} \J{x}^{-1}
\right\| \leq C
\end{align} 
and 
for $\beta < \beta ^{\ast}$ 
\begin{align} \label{aad10}
\left\| 
P  \J{x}^{-1}
\right\| \leq C.
\end{align} 
Here, we only consider the case of $\beta \geq \beta ^{\ast}$. A simple calculation shows 
\begin{align*}
& \J{H}^{2 \beta} x \cdot p (1+p^2+b^2)^{- \gamma} = \J{H}^{2 \beta} (1+p^2+b^2)^{- \gamma} ( (\mbox{bounded opetator}) + p \cdot x) \\ &= 
(\mbox{bounded opetator}) + \J{H}^{2 \beta} (1+p^2+b^2)^{- \gamma}p \cdot x. 
\end{align*}
Moreover, by $[p \cdot x,  \J{p}^{-4 \beta \gamma -1 + 2 \gamma} ] = (\mbox{bounded operator})$, we have 
\begin{align*}
\left\| P\J{p}^{-4 \beta \gamma -1 + 2 \gamma} \J{x}^{-1} \right\|  = 
 C +  \left\| \J{H}^{2 \beta} (1+p^2+b^2)^{- \gamma} \J{p}^{-4 \beta \gamma -1 + 2 \gamma+1}  \right\| \leq C.
\end{align*} 
In the same way, for $1/2 < s \leq 1$, 
$$ \left\| 
\J{P}^s  \J{p}^{-4 \beta \gamma s -s + 2 \gamma s}\J{x}^{-s}
\right\| \leq C. 
$$ 
Now we prove that for $\theta = -4 \beta \gamma s + 2 \beta \gamma -s + 2 \gamma s  \leq -4 \beta s + 2 \beta +s \leq 1$, 
\begin{align} \label{aad11}
\left\| 
\J{x}^{-s} \J{p}^{\theta} \J{H}^{- \beta} \J{x} ^{s} \J{p} ^{4 \beta \gamma s + s -2 \gamma s}
\right\| \leq C
\end{align}
holds if $\beta \geq \beta ^{\ast}$ and for $\theta = 2 \beta \gamma$, 
\begin{align} \label{aad14}
\left\| 
\J{x}^{-s} \J{p}^{\theta} \J{H}^{- \beta} \J{x} ^{s} \right\| \leq C
\end{align}
holds if $\beta \leq \beta ^{\ast}$. Again, we only consider the case of $\beta \geq \beta ^{\ast}$. By Helffer--Sj\"{o}strand's formula (see \cite{HS}), we get 
\begin{align*}
\J{p}^{\theta} \J{H}^{- \beta} \J{x}^s \J{p} ^{4 \beta \gamma s + s -2 \gamma s} = 
\J{p}^{\theta } \J{x}^s \J{H}^{- \beta} \J{p} ^{4 \beta \gamma s + s -2 \gamma s}+ 
\J{p}^{\theta} \CAL{L} \J{p} ^{4 \beta \gamma s + s -2 \gamma s}
\end{align*}
with 
\begin{align*}
\left\| 
\J{p}^{\theta} \CAL{L} \J{p} ^{4 \beta \gamma s + s -2 \gamma s}
\right\| \leq C \left\| 
\J{p}^{\theta -2 \gamma}
\right\| \left\| 
 \J{p} ^{4 \beta \gamma s + s -2 \gamma s - 2 \gamma -1}
\right\| \leq C, 
\end{align*}
where we use $\beta \geq \beta ^{\ast}$, $\theta - 2 \gamma  \leq 2 \beta \gamma -2 \gamma \leq 0 $ and 
$4 \beta \gamma s + s -2 \gamma s -1 \leq s-1 \leq 0$. Moreover, a symbol calculation of the pseudo-differential operator gives 
\begin{align*}
\left\| 
\J{x}^{-s} \J{p}^{\theta} \J{x}^s \J{p}^{-s}
\right\|  \leq C, 
\end{align*}
which proves \eqref{aad11}. By \eqref{aad9} and \eqref{aad11} or \eqref{aad10} and \eqref{aad14}, we get 
\begin{align*}
 \left\| 
\J{x}^{-s} \J{p} ^{\theta} \tilde{\varphi }(H)e^{-itH} \phi
\right\| 
 &\leq 
\left\| 
\J{x}^{-s} \J{p} ^{\theta} \J{H}^{-\beta} \J{P}^{s}
\right\| \left\| 
\J{P}^{-s} \J{H}^{ \beta}\tilde{\varphi }(H)e^{-itH} \phi
\right\|  \\ & \leq C \left\| 
\J{P}^{-s} \J{H}^{ \beta}\tilde{\varphi }(H)e^{-itH} \phi
\right\|
\end{align*}
and, together with \eqref{aad15}, we obtain \eqref{aad16}. 
}
By $\max_{0 \leq \beta \leq 1/2} \theta = \gamma -1/2$, we finally obtain the following smoothing estimate: 
\begin{Thm} 
For the operator $H= a (p^2 +b^2)^{\gamma} +V$ in \eqref{aad12}, 
\begin{align*}
&\sup_{|\lambda| \geq 3R }\left\| 
\J{x}^{-s} \J{p}^{\gamma -1/2}  (H- \lambda \mp 0)^{-1 } \J{p}^{\gamma -1/2} \J{x}^{-s}  \phi 
\right\| \leq C \| \phi \|, \\ 
& \int_{\bf R} \left\| 
\J{x}^{-s} \J{p} ^{\gamma -1/2} \tilde{\varphi} (H)e^{-itH} \phi 
\right\| ^2 dt \leq C \| \phi \|^2.
\end{align*} 
holds. 
\end{Thm} 
If $\gamma =1$, this estimate corresponds to the smoothing estimate for the Schr\"{o}dinger operator with a high-energy cut-off (see, e.g., \cite{EGS}.

\end{document}